\documentclass[letterpaper]{article} 
\usepackage{aaai23}  
\usepackage{times}  
\usepackage{helvet}  
\usepackage{courier}  
\usepackage[hyphens]{url}  
\usepackage{graphicx} 
\usepackage{natbib}  
\usepackage{caption} 
\usepackage{algorithm}
\usepackage{algorithmic}

\usepackage{todonotes}

\setuptodonotes{inline}

\usepackage{multirow}

\usepackage{array}

\newcommand{\strikethrought}[1]{%
  \setbox0\hbox{#1}%
  \usebox0\kern-\wd0\raise.2em\hbox{\vline width\wd0 height 1pt}%
}

\usepackage{amsmath}

\usepackage{newfloat}
\usepackage{listings}
\DeclareCaptionStyle{ruled}{labelfont=normalfont,labelsep=colon,strut=off} 
\floatstyle{ruled}
\newfloat{listing}{tb}{lst}{}
\floatname{listing}{Listing}
\title{Scheduling Ground-Based Telescope Observations with Uncertain Nights}
\author {
    Thomas Rahab Lacroix\textsuperscript{\rm 1}, 
    Pierre Lemaire\textsuperscript{\rm 1}, 
    Anne-Marie Lagrange\textsuperscript{\rm 2,3}, 
    Julien Milli\textsuperscript{\rm 3}, 
    Nadia Brauner\textsuperscript{\rm 1}
}
\affiliations {
    \textsuperscript{\rm 1} Univ. Grenoble Alpes, CNRS, Grenoble INP, G-SCOP, 38000 Grenoble, France\\
    \textsuperscript{\rm 2} LIRA, CNRS, Observatoire de Paris, CNRS, Université PSL, 5 Place Jules Janssen, 92190 Meudon, France\\
    \textsuperscript{\rm 3} Univ. Grenoble Alpes, CNRS, IPAG, F-38000 Grenoble, France\\
    \{thomas.rahab-lacroix, pierre.lemaire, nadia.brauner\}@grenoble-inp.fr,\\ anne-marie.lagrange@obspm.fr, julien.milli@univ-grenoble-alpes.fr
}

\usepackage{bibentry}

\begin{document}

\maketitle

\begin{abstract}
The observation of celestial objects is a fundamental activity in astronomy. Ground-based and space telescopes are used to gather electromagnetic radiation from space, allowing astronomers to study a wide range of celestial objects and phenomena, such as stars, planets, galaxies, and black holes. The European Southern Observatory (ESO) charges each night 83 kEUR \cite{milli2019nowcastingturbulenceparanalobservatory}, so the schedules of the telescopes are really important in order to optimize every second. Ground-based telescopes are affected by meteorological conditions, such as clouds, wind, and atmospheric turbulence. Accurate scheduling of observations in the presence of such uncertainties can significantly improve the efficiency of telescopes use and support from automated tools is highly desirable. In this paper, we study a mathematical approach for scheduling ground-based telescope observations under an uncertain number of clear nights due to uncertain weather and atmospheric conditions. The model considers multiple targets, uncertain number of nights, and various observing constraints. We demonstrate the viability and effectiveness of an approach based on stochastic optimization and reactive strategy comparing it against other methods.
\end{abstract}

\section{Introduction}

This paper is part of several collaborative projects between an Operations Research team and astronomers. The objective is to develop a scheduling software for SPHERE \cite[Spectro-Polarimetric High-contrast Exoplanet REsearch instrument,][]{beuzit2019sphere} installed at the ESO/Very Large Telescope (VLT). This software is designed to take atmospheric uncertainties into account. In the longer term, this tool could be easily adapted for scheduling on any telescope.

The current problem can be summarized as follows: A single telescope is available, and astronomers are interested in a list of $S$ observations with various associated constraints. The time available to carry out these observations is limited, with only $m$ nights at their disposal and the objective is to find a best schedule, taking into account
the specific priority of each target. To quantify this priority, astronomers assign a score $w_j$ to each observation $j$, which represents its gain or scientific value. The goal is therefore to select observations and schedule them over $m$ nights, while maximizing the sum of the gains of the completed observations.

In general, only one observation can be done at a time with a single telescope. The schedule is subject to astronomical constraints: each observation requires a specific amount of time to be completed and due to the Earth's rotation, observations are not always feasible; they have visibility windows that vary from night to night. Those astronomical constraints define a complex, yet classic, scheduling problem. In addition meteorological and atmospheric conditions are important constraints for ground-based telescopes conversely to space telescopes. The main contribution and originality of the present paper is to consider these additional constraints. We assume as a first step that when the weather is bad, the telescope has to be shut down; as the weather cannot be known beforehand, it means that the actual number of available nights is stochastic.

In the remainder of this paper, we first synthetize the state-of-the-art of scheduling for astronomical observation and with uncertain number of nights. Then, we formalize our problem and propose solution methods based on constraint programming and stochastic optimization with a discussion on their implementation. Numerical experiments based on artificial instances are presented to study the tractability of the problem, the behavior of the methods and the gains obtained. We conclude with a summary of the main results and a discussion about the next step of the project.

\section{Related Work}

\subsection{Scheduling in an Astronomical Context}

Scheduling observations is an old challenge. Early works handle it with rule based AI \cite{MILLER1987301} for the Hubble Space Telescope (HST) and \cite{Johnston1988AutomatedOS} for the VLT. Then Spike \cite{johnston1990spike} was developed for HST including also methods for ground-based telescopes \cite{johnston1994spike}. Recently STScI, the institute developing the tools used for the James Webb Space Telescope, implemented the improvements made to the HST's scheduler \cite{adler2014planning}.

The current project started in 2015 for the scheduling of 200 nights in a survey dedicated to the search for exoplanets with the SPHERE instrument. \citet{brauner:hal-01237688} introduces a problem with specific constraints on top of the classical astronomical constraints. Then, \citet{paper_catusse} proposes a branch and price approach to solve it, efficient due to the very specific structure of the nights. Next, \citet{Lagrange2016SPOTAO} presents the concrete implementation of this method in a software for astronomers. Finally, the thesis of \citet{these_flo} addresses the problem from a more academic perspective and considers reducing slightly the observing time of targets, which also reduces the gains but ultimately may allow for more observations. None of these works takes weather uncertainty into account.

In the case of scheduling Earth observations from satellites, weather conditions is first considered in \citet{uncertain_model}. The main idea is to compute a new weight of interest for each observation, depending on realization probabilities and remaining observation opportunities. \citet{Liao2005SatelliteIO} is the first practical implementation of satellite observations with uncertainties. However the forecast of the weather is assumed to be good and the long term is not considered. \citet{he_2015} model a new benefit which depends on the cloud coverage. They solve their problem for only one satellite orbit and well known meteorological conditions. In the case of astronomical observations with ground-based telescopes, we must consider several nights and highly uncertain meteorological and atmospheric conditions. Indeed, those conditions directly impact the quality of the astronomical data. Current works aim at modeling the astronomical image quality as a function of these  external conditions. \citet{WANG2021107292} use budgeted uncertainty to describe the cloud coverage. They also improve the solving part with column generation and a simulated annealing. 

As a summary, the scheduling of Earth observations with satellites is similar in many ways to the scheduling of celestial observations with ground-based telescopes but the weather parameters are different: the meteorological conditions are more complex in the case of astronomical observations; in addition, for astronomical observations with ground-based telescopes, atmospheric optical turbulence conditions are also to be taken into account and are fully uncertain. The horizons involved are also different so the scheduling of ground-based astronomical observations requires a specific approach.

\subsection{Scheduling for Throughput Maximisation}

For consecutive nights, the positions of the stars in the sky at a given time do not change much. The visibility windows are only slightly shifted. In such a case, we can assume that the visibility windows remain the same each night, making the nights identical.

In scheduling theory terms, we consider stars as tasks to be accomplished and nights as independent parallel machines. Therefore in Graham's notation, the problem is a $P_m|r_j|\sum w_j U_j$, \textit{i.e.}, we want to schedule tasks on $m$ parallel identical machines with release dates and deadlines and the objective is to maximize the weighted throughput (weighted number of scheduled task).
The problem is NP-hard; indeed even restrictions already are, for instance with a single night and no release dates \cite[$1||\sum w_j U_j$,][]{lawler1969functional}; with no weights $(1|rj,dj|.)$. It implies that even restricted cases can already be computationally challenging.

\subsection{Uncertain Number of Machines}

Some nights, planned observations might not be feasible because of meteorological and astronomical conditions. For simplicity, we consider that, for each night, either the weather is clear (allowing all observations) or it is bad, in which case the telescope is completely shut down until the next night. The weather conditions are moreover assumed to be independent from one night to another (which is also an approximation). Therefore, we consider that the number of available nights $m$ is uncertain.

Scheduling with uncertainty has been well studied as any practical application has some part of uncertainty; however, this problem has mainly been addressed through the uncertainty on the tasks. \citet{10.1145/3340320} study scheduling with an uncertain number of machines. \citet{balkanski2022schedulingspeedpredictions} generalizes it to uncertainty on the speed of the machines. The uncertain number of machines can then be seen as the special case with binary uncertain speeds. Those papers, consider some type of robust optimization (see also \citet{eberle2022speedrobustschedulingsand}).

In our context with repeated use of the scheduler, the worst case approach of robust optimization is not the most appropriate. \citet{buchem2024schedulingstochasticnumbermachines} and \citet{epstein2024efficientapproximationschemesscheduling} maximize the expected value in a context of stochastic optimization.

All articles addressing uncertain number of machines are mainly focused on Polynomial Time Approximation Scheme (PTAS). The objective functions are mainly based on completion times like makespan \cite{minarik2024speedrobustschedulingrevisited}, minimal load, weighted completion time \cite{pmlr-v202-lindermayr23a} and no feasibility windows have been considered.

\section{Problem Statement}

In the scheduling framework, our problem can then be described as the scheduling of tasks with release dates and deadlines and with an uncertain number of machines maximizing the expected value of the tasks scheduled, \textit{i.e.} $P_m | r_j | \sum w_j U_j$ with uncertainty on the number of machines $m$.

\subsection{Problem Formulation}

The problem is defined by the following data:

\begin{itemize}
    \item Total number of nights: $M$
    \item The probability that exactly $m$ nights are observable: $\pi_m, m \in [\![0,M]\!]$
    \item Number of observations: $S$
\end{itemize}

\noindent Independently of the night, for each observation $j \in [\![1,S]\!]$:

\begin{itemize}
    \item Release date: $r_j$
    \item Deadline: $d_j$
    \item Processing time: $p_j$
    \item Gain: $w_j$
\end{itemize}

Surprisingly, $m$ is not really a data of the problem. $M$ is the total number of nights and is just here to bound $m$. In fact $m$ is just a tacit consequence of the probabilities $\pi_m$.

It is important to note that the data for $\pi_m$ represents the probability of having exactly $m$ nights available, not the probability of a specific night $m$ being observable. For example $\pi = (0.1, 0.3, 0.2, 0.2, 0.2)$ means that there is a probability of 0.1 that exactly zero night is observable, a probability of 0.3 that exactly one night is observable and so on. While the detailed probabilities for each night could allow the recovery of $\pi_m$, the reverse is not true. By working with fewer details, the model focuses on minimal but sufficient information, which is an advantage.

Given the problem data, the probability distribution, and the repetitive nature of this algorithm's application, modeling it through stochastic optimization to maximize the expected total gain is appropriate.

The objective of the model is to maximize the Expected Total Gain (ETG):

\begin{equation}
  \max \displaystyle \sum_{m=1}^{M} \pi_m \cdot \displaystyle \sum_{i=1}^{m} \displaystyle \sum_{j=1}^{S} w_j \cdot b_{ij}
  \label{eq1}
\end{equation}

\noindent where $b_{ij}$ is 1 if observation $j$ is scheduled during night $i$, and 0 otherwise. All other constraints are classical scheduling constraints: observations do not overlap within one night, and each observation is scheduled in at most one night, lasts $p_j$ units of time, and is carried out within its release date and deadline.

\subsection{Dominance Analysis}

Since the specific night on which observations are scheduled does not matter, solutions equal, up to the order of the nights, are equivalent. Indeed, a set of solutions with exactly the same scheduling but with swapped nights is dominated by some of its members: given the uncertainty surrounding distant nights, it is better to prioritize nights that yield the highest gains as early as possible.

Consequently, all solutions that do not follow a decreasing night gain order can be eliminated, as they are dominated. This approach significantly reduces the solution space without losing any relevant solutions.

This property underlines the fact that what matters is not the absolute position of a night but its relative importance compared to the other nights. It aligns well with the modeling choice that $\pi_m$ is the probability of having exactly $m$ nights, and not the probability of a specific night $m$ being observable.

\subsection{Significance of Uncertainty}

\begin{figure}[h!]
\centering
\includegraphics[width=0.5\textwidth]{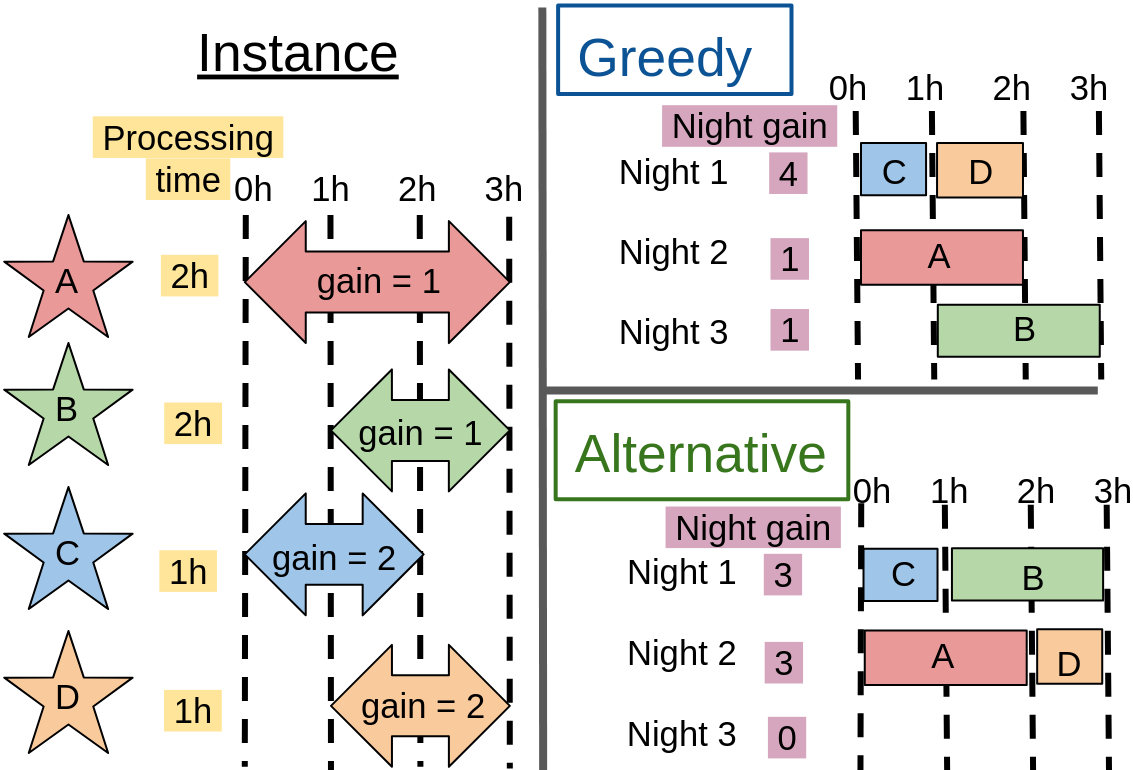}
\caption{\label{fig:ctr_ex_greedy_opti}Counterexample data and solutions with 4 stars.}
\end{figure}

Even though sorting nights by decreasing gain is dominant, a greedy algorithm may fail to provide an optimal schedule. Indeed, Table~\ref{tab:signi_uncert} and Figure~\ref{fig:ctr_ex_greedy_opti} describe a counterexample involving 3 nights and 4 stars. With all three nights available, both solutions are optimal since all the stars are observed. However, if only two nights are observable, the Alternative solution achieves a total gain of 6, strictly greater than the Greedy solution’s gain of 5. On the other hand, if the weather is very poor and only one night is observable, the Greedy solution is better, with a gain of 4 compared to the Alternative solution's gain of 3.

This example also demonstrates that the probability distribution over the number of observable nights has a significant impact on the optimal solution. 

\begin{table}[h!]
    \centering
    \begin{tabular}{c|c|c|c|c}
        Star & $r_j$ & $d_j$ & $p_j$ & $w_j$ \\\hline
        A & 0 & 3 & 2 & 1 \\
        B & 1 & 3 & 2 & 1 \\
        C & 0 & 2 & 1 & 2 \\
        D & 1 & 3 & 1 & 2
    \end{tabular}
    \caption{\label{tab:signi_uncert}Counterexample data with 4 stars, release dates ($r_j$), deadlines ($d_j$), processing times ($p_j$) and gain ($w_j$).}
\end{table}

\section{Solution Methods}

Three algorithms where developed on this problem to compare their solutions. One which considers uncertainty, a greedy and an oracle.

\subsection{Stochastic Optimization}

The stochastic approach seeks to maximize the expected value of the schedule as described in expression (\ref{eq1}). It is based on a mathematical formulation of the problem (eg., Constraint Programming or Mixed Integer Linear Programming) and it is solved using an adequate solver. This method can yield exact solutions at the expense of high computation times.

\subsection{Greedy Heuristic}

To our knowledge, this approach is the most widely used strategy for astronomical scheduling. Each night is planned as if it were the last, prioritizing the best remaining observations.

While simple and very fast, this heuristic can easily fail, as illustrated in Figure~\ref{fig:ctr_ex_greedy_opti}. A key question is thus to test the gain brought by the Stochastic Optimization over the Greedy Heuristic, with regard to the extra computing effort required.

\subsection{Omniscient Approach}

Unlike the other two approaches, the Omniscient algorithm does not provide an actual solution. It assumes perfect knowledge of the future, removing any uncertainty. So it computes the perfect schedule for each number of nights. It corresponds to what an Oracle could do, and it is used to measure the cost of uncertainty.

\section{Implementation}

The implementation is based on Constraint Programming because it is efficient to manage no-overlap constraints. The Choco constraint programming solver is used, because it is efficient, widely used, actively maintained and a free software.

The implementation in the Java Choco module of the Stochastic Optimization is straightforward, by just modeling the scheduling constraints using, for instance, the no-overlap constraint together with the objective function as in expression~(\ref{eq1}). The Omniscient Approach is calculated by aggregating the results over $M$ executions, where in each of the $M$ configurations exactly one of the $\pi_m$ is set to 1. Similarly, for the Greedy Heuristic, the instances are modified by setting the number of nights to 1, setting the probabilities to 1 for that single night, then removing the scheduled observations from the instances and repeating this process $M$ times.

The implementation has not been optimized for performance, as our primary focus is on behavior and order of magnitude.

\subsection{Instance Generation}

Instances are essentially defined and randomly generated based on four meta-parameters. These four meta-parameters are $M$, $S$, $lenNight$ and $maxGain$. $M$ and $S$ directly correspond to the problem data, representing the number of nights and the number of observations, respectively. $lenNight$ represents the length of the nights. For each star, the release date and deadline are drawn uniformly among integers between 0 and $lenNight$, ensuring that the release date is strictly less than the deadline. The processing time is then randomly drawn among integers between 1 and deadline minus release date. $maxGain$ determines the gain associated with each star, which is randomly drawn among integers between 1 and $maxGain$.

Probabilities for each number of available nights are generated as follows. An initial value for $m$ nights ($q_m$) is randomly drawn between 0 and 1. The actual probability ($\pi_m$) is then calculated as: $\pi_m = \frac{q_m}{\displaystyle\sum_i q_i}$. Finally, the probability for $\pi_0$ is adjusted to correct rounding errors, ensuring that the total probabilities sum to exactly 1. Since $\pi_0$ has no impact on the instance or its optimal solution, modifying its value is not problematic.

The advantage of this method is that it allows the generation of any possible instance. However, these instances are not generated with equal probabilities.

\subsection{Redundant Constraints to Speed Up the Resolution}

In order to speed up the resolution through constraint programming, we propose the three following redundant sets of constraints.

\begin{enumerate}
    \item Decreasing Gain (DG). This constraint ensures that the gain decreases night after night. Indeed, if a solution does not satisfy this constraint, swapping the two nights in question would strictly improve the overall objective-function. This constraint implements the concept of dominance detailed above. Formally, it states that : \\ $NightGain(i-1) \geq NightGain(i), \forall i \in [\![2,M]\!]$.
    
    \item Bounded Observations (BO). $boundMaxObs$ is an upper bound on the number of observations possible in a single night. This constraint simply ensures that the number of observations does not exceed the calculated bound. Formally, it states that : \\ $\displaystyle\sum_{j=1}^{S} b_{ij} \leq boundMaxObs, \forall i \in [\![1,M]\!]$ \\ There are several ways to determine this bound. It is necessarily less than or equal to $lenNight$ (achievable if enough observations have a processing time of 1). Here, $boundMaxObs$ is calculated as the number of observations with the shortest processing times, whose cumulative sum is less than or equal to $lenNight$. In practice, however, it is not $lenNight$ that is used, but rather the latest due date minus the earliest release date among all observations. Aside from this detail, observation windows are not considered at all. If they were included, the bound could be even stricter, but this would involve solving a problem that is NP-complete during the presolve phase ($1|r_j|\sum U_j$). 
    
    \item Increasing Cumulative Gain (ICG). This constraint simply ensures that the scheduled nights are worth doing setting that the cumulative gain is increasing. By construction, any solution satisfy this property, but stating it explicitly may guide efficiently the solver. Formally, it states that : \\ $CumulGain(i-1) \leq CumulGain(i), \forall i \in [\![1,M]\!]$
\end{enumerate}

The improvement in computation time thanks to these three constraints is detailed in Table~\ref{tab:cut_eval}. It presents the mean resolution times for 100 instances with 4 nights, 20 observations, a $lenNight$ of 5, and a $maxGain$ of 10, for all configurations (whether the constraints are applied or not). The size of the instances was chosen to have non trivial solutions and large but acceptable computation times.

\begin{table}[h!]
\centering
\begin{tabular}{c|c|c|r r r}
 & & & \multicolumn{3}{c}{Solving time (s)} \\  
DG & BO & ICG & Mean & [ Min , & Max ]   \\\hline
$\times$ & $\times$ & $\times$ & 38.9 & 0.4 & 407.3 \\
$\times$ & $\times$ & & 39.1 & 0.4 & 454.1 \\\hline
$\times$ & & $\times$ & 181.9 & 0.7 & 5001.2 \\
$\times$ & & & 183.3 & 0.7 & 4954.4 \\\hline
 & $\times$ & & 928.0 & 1.4 & 16392.7 \\
 & $\times$ & $\times$ & 978.7 & 1.4 & 18239.6 \\
 & & & 1131.7 & 1.2 & 12970.8 \\
 & & $\times$ & 1148.1 & 1.1 & 16224.1
\end{tabular}
\caption{\label{tab:cut_eval}Redundant constraints evaluation. Decreasing Gain (DG), Bounded Observations (BO) and Increasing Cumulative Gain (ICG).}
\end{table}

It is clear that the Decreasing Gain configuration significantly accelerates computation times. It is undeniably useful and effectively divides the table into two groups: configurations with mean computation times over 900 seconds and those under 200 seconds.

Next, Bounded Observations appears to be genuinely useful. It combines very well with Decreasing Gain to separate the configurations with mean computation time over 180s and those under 40s.

Finally, the least useful constraint seems to be Increasing Cumulative Gain which is worst alone than all other configurations including the one with no redundant constraint. However, with one or both of the other constraints, the mean computation is slightly reduced by Increasing Cumulative Gain. The last differences being very small, it does not allow to conclude that this constraint is effective but it does not seem counterproductive. Therefore, it is retained for the remainder of the experiments.

\section{Experiments and Results}

\subsection{Size of Solvable Instances}

The next step is to generate instances that will cover the feasible solution space as much as possible and solve them to have an idea of the capacity of the solver.

The first experiment is to challenge the solver capability by running it on instances from trivial to more complex. To this extent, using the methodology explained above, we generated instances with the following meta-parameters: 
\begin{itemize}
    \item $M$: 1, 2, 5, 10, 20
    \item $S$: 1, 2, 5, 10, 20, 50, 100
    \item $lenNight$: 1, 2, 5, 10
    \item $maxGain$: 1, 2, 5, 10
\end{itemize}

For each of the 560 unique meta-parameter categories, 10 instances were generated to slightly explore random variations within each category for a total of 5,600 instances.

For each of them, a time of 10 seconds is allocated for the resolution before moving to the next instance. This time constraint provides a quick overview of the algorithm's computational capacity to determine the upper limit of solvable instance sizes. The reasons of those 10 seconds is that the total computation time is 15 hours if all instances require 10 seconds each. For future experiments, when there is less instances, more computation time is allocated in order to reproduce the real operational conditions.

Here are some interesting statistics obtained after the resolution. They are 3922 instances among the 5600 (70.0\%) solved to optimality. The others would need more than 10s to finish their computation. Then, only 84 instances over the 3922 (2.1\%) have a gap between the ETG (Expected Total Gain) obtained with the Stochastic Optimization and the Greedy Heuristic. This value is low because many instances have a low number of observations compared to the number of nights so they are all scheduled easily even with the Greedy Heuristic. Finally, among the 84 instances the mean upgrade of Stochastic over Greedy is 4.9\% --- where the upgrade of A1 over A2 is defined by:

\begin{equation}
Upgrade(A1,A2) = \frac{ETG[A1]-ETG[A2]}{ETG[A2]}
\label{eq2}
\end{equation}

It means that the Greedy Heuristic is beaten by the Stochastic Optimization by 4.9\% of its ETG. Because of the various instances sizes, this value is not really pertinent and is just here to give an idea of the possibilities. This upgrade varies from almost 0\% to 66.4\% on some instances.

After testing the solver on all type of instance sizes, it is time to find a size that provides instances that tend to be similar to real astronomical instances. A $lenNights$ of 5 is close to the number of observations done in one night. Same for $maxGain$ of 10 which is close to the scale used currently by the astronomers. Then to keep the problem non trivial, enough observations must be in the pool and having $S \geq M \times lenNights$ allows this. So after few more tests and in order to have an acceptable computation time (around few minutes) for one instance, the appropriate solution size selected is 4 nights, 20 observations, 5 $lenNights$ and 10 $maxGain$. So most of the following experiments will be on instances with a size close to this.

\subsection{Impact of Probabilities on Computation Times}

Another factor affecting computation time is the probabilities vector. Table~\ref{tab:proba_affect_comput_time} shows the mean, min, and max solving times for 100 instances with 4 nights, 20 observations, 5 $lenNight$, and 10 $maxGain$. Those 100 instances are each solved 4 times (one for each line of Table~\ref{tab:proba_affect_comput_time}) with their probabilities modified to correspond to the most extreme situations. As expected, the higher the probability of an instance having a large number of observable nights, the more difficult the instance becomes.

\begin{table}[h!]
    \centering
    \begin{tabular}{c|r r r}
     & \multicolumn{3}{c}{Solving time (s)} \\  
    Probabilities $\pi$ & Mean & [ Min , & Max ]   \\\hline
    (0,1,0,0,0) & 0.005 & 0.001 & 0.054   \\
    (0,0,1,0,0) & 0.180 & 0.015 & 1.270  \\
    (0,0,0,1,0) & 9.058 & 0.240 & 126.945 \\
    (0,0,0,0,1) & 432.664 & 3.124 & 7057.913 
    \end{tabular}
    \caption{\label{tab:proba_affect_comput_time}Link between probabilities and solving time.}
\end{table}

Computation times seem to increase exponentially as high probabilities are assigned to a large number of observable nights. For a real instance where probabilities are not binary, it is difficult to estimate computation time solely based on probabilities. However, some tests provided a rough idea. The computation time for an instance with randomly generated probabilities will be approximately 10 times lower than the computation time on the same instance with probability set to 1 for the largest number of nights.

\subsection{Behavior of the Algorithms}

Now that we have the various sizes of solvable solutions, the next step is to examine what the solutions look like. For the same instance, we will compare the Greedy Heuristic with the Stochastic Optimization. We will also compare the Stochastic Optimization with the Omniscient Approach, which is impossible to achieve, and simply shows the maximum feasible with $m$ nights. We focus on the total gain and its progression based on the number of nights.

Since the last step of the Omniscient Approach corresponds to the last raw of Table~\ref{tab:proba_affect_comput_time}, more computation time will be required to solve each instance properly. As the Greedy Heuristic is very fast, its computation time is negligible compared to the other algorithms. In the end, the total computation time to solve properly 100 instances each having the size 4 nights, 20 observations, 5 $lenNights$ and 10 $maxGain$ with all three algorithms is 13h20.

\begin{figure}[h!]
\centering
\includegraphics[width=0.5\textwidth]{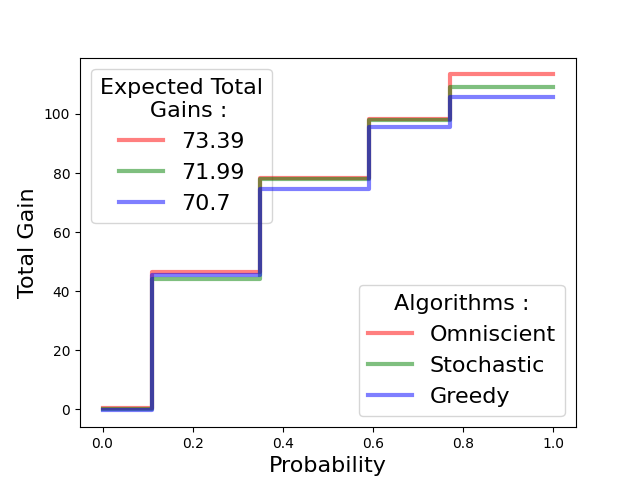}
\caption{\label{fig:aire_sous_courbe_style}Evolution of the total gain for the three algorithms on an instance with 4 nights.}
\end{figure}

Figure~\ref{fig:aire_sous_courbe_style} shows a typical pattern of the algorithms behaviour on an instance.

The vertical axis of this figure represents the cumulative gain over the nights. On the horizontal axis, we have the probability of observable nights. Note that each step corresponds to one additional favorable night. This representation was chosen because the area under a curve equals the objective function. Indeed, the height of each step is the total gain achieved if the corresponding number of nights is observed, and the width of each step is the probability of that happening. Thus, the area under a curve precisely represents the average total gain achieved by the solution for the corresponding instance. The curves in Figure~\ref{fig:aire_sous_courbe_style} can also be seen as the inverses of the cumulative distribution functions of the total gain.

Here are some important remarks before the detailed analysis of the curve. The Omniscient Approach curve might be impossible to match: all other curves will necessarily be below or equal to it. Indeed, the Omniscient Approach does not represent a single schedule but the best possibility among all schedules for each step/number of nights. There is no guarantee that the transition between its steps can be achieved by a single schedule.

The Greedy Heuristic is unbeatable for a single night. Its curve will perfectly overlap with the Omniscient Approach curve for at least the first night of all instances. By definition, the Greedy Heuristic maximizes its gain for the current night, then maximizes again for the following night, and so on. If the Stochastic Optimization curve is below the Greedy Heuristic during the initial nights, this is normal. Stochastic Optimization groups observations that fit well together and reserves them for later. It lags slightly at the beginning but catches up as the Greedy Heuristic slows down in terms of well-filled nights, and Stochastic Optimization continues to improve and surpasses the Greedy Heuristic. Ultimately, what matters is the area under the curve. Even if Stochastic Optimization seems to make sacrifices at the start, in reality, it  wins in the long run.

The gap between the Omniscient Approach curve and the Stochastic Optimization curve perfectly represents the cost of uncertainty. If there were no uncertainty, one step of the Omniscient Approach curve would be realized because it is optimal. However, due to this uncertainty, the Stochastic Optimization curve, which maximizes the expected gain, can be considered optimal for the problem with $m$ uncertain. The gap between the two curves thus perfectly quantifies how much uncertainty causes us to lose, \textit{i.e.}, the cost of uncertainty.

\subsection{Quantitative Analysis of the Algorithms}

Here are some interesting statistics obtained after the resolution. They are 43 instances among the 100 with a gap between the ETG (Expected Total Gain) obtained with the Stochastic Optimization and the Greedy Heuristic. Moreover, among the 43 instances, the mean improvement of Stochastic over Greedy is 62\% with the improvement of A1 over A2 being:

\begin{equation}
Impro(A1,A2)=\frac{ETG[A1]-ETG[A2]}{ETG[Omniscient]-ETG[A2]}
\label{eq3}
\end{equation}

It means that the Stochastic Optimization achieves 62\% of the possible improvement. It closes 62\% of the gap between the Greedy Heuristic and the Omniscient Approach.
The last percents remaining are due to the uncertainty. Using this metric for 5600 instances is impossible because many of those instances are unsolvable with the Omniscient Approach. This improvement is great and shows the interest of working on such methods. Moreover, this improvement is expected to increase on larger instances.

Other interesting statistics are the ones linked to the cost of uncertainty. There are 56 instances among the 100 with a gap between the ETG obtained with the Omniscient Approach and the Stochastic Optimization. Among those 56 instances, the mean upgrade of Omniscient over Stochastic (see expression (\ref{eq2})) is 0.9\%.

It means that the remaining upgrade possible from the Stochastic Optimization to the Omniscient Approach is 0.9\% of the Stochastic Optimization. So the cost of uncertainty is here quite low. Unfortunately, this upgrade is expected to increase on larger instances.

Those results are very encouraging. These curves are successful because we observe the Greedy Heuristic globally below the Stochastic Optimization curve, and the Omniscient Approach curve is always above. This demonstrates the importance of approaching the problem with more complex algorithms than the current Greedy  Heuristic.

\section{Reactive Strategy}

In the astronomy problem, the identical and independent machines represent nights. They are spread over time, and the problem's data fluctuate. Thus, it seems pertinent to adjust schedules for these two reasons: there is time to recalculate between two nights, and changes in the data after execution could lead to different optimal decisions.

So in order to do so, a Reactive Strategy calculates the optimal stochastic schedule, after each night, removing the observations that where effectively done until then. It is a rolling horizon method \cite{lv2023autonomous} that benefits from an update of the real data after each night/iteration. After measuring the cost of uncertainty, which is visualized by the gap between the Omniscient Approach and the Stochastic Optimization, we now focus on measuring the gap between a static solution and a reactive solution, that corresponds to the cost of stubbornness.

\subsection{Evolution of the Data in the Reactive Context}

For our experiments, a complex aspect to consider is the evolution of the data. Throughout a rolling horizon, the subsequent state must be consistent with the previous one to remain realistic. We assume here that the data for observations (release date, due date, processing time, and gain) do not change. 

After the execution of each night, the number of observations will evolve according to what has been completed, and the number of nights will decrease by one. The last data point is probability. Previously, probabilities were drawn completely at random, but maintaining consistency between the old probability list and a new one with one fewer element is difficult.
The chosen solution to address this problem is to assume that probabilities follow a binomial distribution with the two parameters being the remaining number of nights $m$ and the probability for one night to be observable. This approach is fairly realistic since the independent nights all share the same likelihood of being observable. For the  data, we retain a system that provides the probability of observing $m$ nights, which will be recalculated at each step using the binomial distribution.

\subsection{Reactive Scenario Simulation}

To experiment with this idea, a simulator was developed. It closely looks like the previously developed algorithms since it operates in the same way. The only difference is that it recursively explores all possible scenarios for the realization of nights, resolving the instance at each step, planning the first night, removing the completed observations from its data, and continuing the exploration. Each scenario is then resolved using the rolling horizon method.

Since each night is either observable or not, this boils down to exploring a binary tree where the depth corresponds to the number of nights. Therefore, there are $2^M$ different scenarios, which correspond to the leaves of the tree. The solver will be used $2^M-1$ times (there is no need to reuse the solver if the current night is not observable).

Regarding computation times, the first branches of the tree clearly require more time than the last ones. As previously observed, computation time increases exponentially with the number of nights. Thus, the very first solver call accounts for approximately 90\% of the total computation time. It is worth noting that this initial solution is exactly the Stochastic Optimization solution for the non-reactive problem. Therefore, even though this has no practical significance because the executions are separated by 24 hours, from an experimental perspective, the total computation time increases only slightly between solving the non-reactive and reactive problems.

\subsection{Reactive Strategy Evaluation}

In the reactive strategy simulation, for a given instance, each leaf of the tree (\textit{i.e.} each scenario) has a total gain and a probability to occur thus allowing to calculate the expectation of the Reactive Strategy for this instance. This value is compared to the expectation of the Stochastic Optimization for the same instance.

The instances chosen for this test campaign are as follows: $M$ = 4, $S$ = 20, $lenNight$ = 5, $maxGain$ = 10 with 100 distinct instances. For each instance, the probability for a night to be observable is drawn uniformly between 0 and 1. 

Table~\ref{tab:study_100_instances} presents many interesting statistics obtained after the resolution of the 100 instances.

\def\myraisebox#1#2{\leavevmode\vbox to 0cm{\kern-4mm\hbox to 0cm{\hss\hbox{#2}\hss}\vss}}
\begin{table}[h!]
    \centering
    \begin{tabular}{l|c|c}
    & Proportion with  & Mean  \\  
    & a not null gap & improvement \\
    \hline
    A1 = Stochastic & & \\
    A2 = Greedy & \myraisebox{2mm}{35\%} & \myraisebox{2mm}{75\%}  \\
    \hline
    A1 = Omniscient & & 100\% \\
    A2 = Stochastic & \myraisebox{2mm}{45\%} & (upgrade 0.7\% )\\
    \hline
    A1 = Reactive & & \\
    A2 = Stochastic & \myraisebox{2mm}{58\%}  & \myraisebox{2mm}{57\%}  \\
    \hline
    A1 = Omniscient & & 100\% \\
    A2 = Reactive & \myraisebox{2mm}{69\%}  & (upgrade 0.4\%)
    \end{tabular}
    \caption{\label{tab:study_100_instances}
    Statistics on 100 instances. The gap is defined as $ETG[A1]-ETG[A2]$. The improvement is defined in expression (\ref{eq3}). The upgrade is defined in expression (\ref{eq2}).}
\end{table}

The first two rows of Table~\ref{tab:study_100_instances} do not concern the Reactive Strategy and show the exact same statistics as the one presented in the static version of the problem. The only thing that changed between those values are the probabilities generated for the instances. They were either drawn uniformly either following a binomial distribution. It appears that in Table~\ref{tab:study_100_instances} there are less instances with a not null gap but they have larger mean improvement or upgrade.

The last two rows of Table~\ref{tab:study_100_instances} allow to see the performances of the Reactive Strategy. By being better than the Stochastic Optimization on 58 instances among the 100 and in those 58 instances being able to close the remaining gap between the Omniscient Approach and the Stochastic Optimization by 57\%, the Reactive Strategy demonstrates its effectiveness. However, there are still 69 instances over the 100 with a gap between the Omniscient Approach and the Reactive Strategy but the mean upgrade of 0.4\% on those 69 instances is very low, proving that the uncertainty is almost fully managed.

\begin{figure}[h!]
\centering
\includegraphics[width=0.5\textwidth]{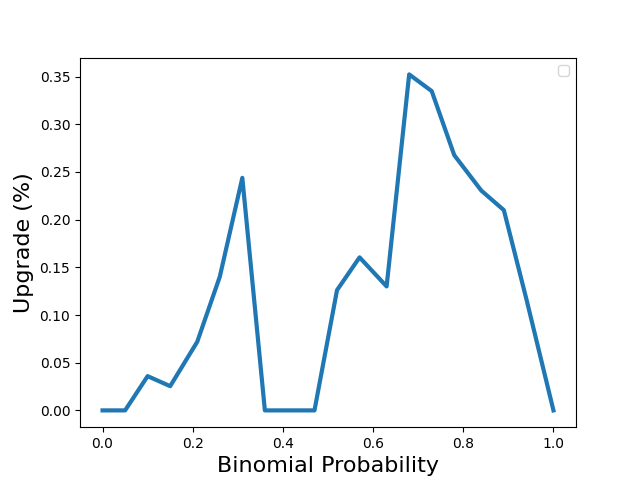}
\caption{\label{fig:evolu_esp_proba}Evolution of the upgrade of Reactive over Stochastic (see expression (\ref{eq2})) depending on the binomial probability on one instance with $M = 4$, $S = 20$, $lenNight = 5$, $maxGain = 10$.}
\end{figure}

Finally, the last test conducted in this section on Reactive Strategy is a study of the impact of the success probability of the binomial distribution. Depending on the value of this parameter, what will the upgrade of Reactive over Stochastic be? Thus, for only one of the 100 instances, this exact test was conducted. The result can be seen in Figure~\ref{fig:evolu_esp_proba}. When the success probability of the binomial distribution is at 0 or 1, then naturally, the Stochastic and Reactive solutions will be the same, therefore, it is normal for the upgrade to be exactly 0\%. It is interesting to note two bumps on the curve. When the binomial probability is close to 0.5 there is no upgrade and there are two pics respectively around 0.25 and 0.75 of binomial probability. It can be interpreted as when an unexpected scenario with the Stochastic Optimization is happening (the two pics), the Reactive Strategy is able to correct it and diverges from the Stochastic Optimization. However, when the binomial probability is close to 0.5, all the scenario are expected and whatever the realisation, the Reactive Strategy is not able to create a separation with the Stochastic Optimization.

\section{Conclusion}

We consider the problem of scheduling astronomical observations from a ground-based telescope taking into account the external (meteorological and atmospheric turbulence) conditions. This is modeled as a classical parallel machine scheduling problem with time windows and the weighted  throughput as the objective, the observations being the tasks and the nights being the machines. The specificity of this work is to consider uncertain weather conditions, modelled as an uncertain number of machines/nights. The approach is based on stochastic optimisation using a constraint programming solver. It is compared to a greedy approach and to the optimal approach when uncertainty is revealed. 

This study has demonstrated the benefits of taking explicitly uncertainties due to weather conditions into account. Reactive strategies that re-optimize stochastic schedules are particularly promising. Besides, trends are such that gains can be expected to be even larger for bigger and more realistic instances.

Beyond these findings, this work serves as a necessary preliminary step, helping to identify the key aspects that should be prioritized in future research. By building on these insights, the next phase of work can focus on refining the methodology and tackling the remaining challenges to further enhance both theoretical and practical outcomes.

\section{Future Work}

From a theoretical perspective, transitioning from Constraint Programming (CP) and Choco to Mixed-Integer Linear Programming (MILP) combined with state-of-the-art solvers would enable handling larger instances, pushing the limits of practical solvability. Furthermore, that would allow to consider non identical nights. Additionally, refining our calculations and visual representations should help highlight key trends more distinctly.

Another important step is to relax some of the simplifying assumptions made in this study. In particular, incorporating mid-term weather forecasts and short-term atmospheric turbulence predictions would provide a more realistic framework. Instead of a binary classification of meteorological conditions, a quantitative approach (wind speed and direction, among other variables) would enhance the model's applicability.

Moreover, increasing the problem’s complexity by considering five categories rather than just two (e.g., clear, poor conditions, southern wind, northern wind, and moderate conditions) would better capture the diversity of real-world scenarios. A step even further will be to refine, for each observation, the image quality depending on the forecast conditions. Another possible improvement would be to update the observation list during the rolling horizon; for example by re-observing certain targets when they prove to be very interesting.

Finally, applying the methods developed here to real astrophysical data would serve as a critical validation step, bridging the gap between theory and practice.

By pursuing these directions, we aim to refine both the theoretical and empirical understanding of the problem, ultimately leading to more scalable solutions and better fitted to the needs of the astronomers. 

\newpage

\section{Acknowledgements}

This work is part of the PEPR ORIGIN ("Smart Scheduling" work package)  supported by the {\em Agence Nationale de la Recherche} under the France 2030 program (ANR-22-EXOR-0018).

\bibliography{aaai23}

\end{document}